\documentclass[11pt]{article}
\edef\qedrestoreat{\noexpand\catcode\lq\noexpand\@=\the\catcode\lq\@}
\catcode\lq\@=11

\ifx\protect\undefined\let\protect\relax\fi

\def\qed{\protect\@qed{$\qedsymbol$}}
\def\pushright{\protect\@pushright}

\def\QED{\protect\@qed{{\rm Q.E.D.}}}

\def\QEI{\protect\@qed{{\rm Q.E.I.}}}

\def\Proof{\protect\@Proof}\def\endProof{\protect\@endProof}%
\let\proof\Proof
\def\Proofof#1{\protect\@Proofof{#1}}\def\endProofof{\protect\@endProofof}%
\let\proofof\Proofof\let\endproofof\endProofof

\def\qedsymbol{\raisebox{-.2ex}{$\Box$}}

\def\TheWordProof{\em Proof.}
\def\TheWordProofof#1{\em Proof of #1.}

\def\ProofFont{}

\newif\ifAutoQED\AutoQEDfalse

\newif\ifNumberResults
\ifx\theorem@style\undefined
   \NumberResultstrue
   
\fi

\def\parag@pushright#1{{
    \parfillskip=0pt            
    \widowpenalty=10000         
    \displaywidowpenalty=10000  
    \finalhyphendemerits=0      
    \hbox@pushright             
    #1
    \par}}

\def\hbox@pushright{
    \unskip                     
    \nobreak                    
    \hfil                       
    \penalty50                  
    \hskip.2em                  
    \null                       
    \hfill                      
}%

\newif\if@qed\@qedfalse
\def\save@set@qed{\let\saved@ifqed\if@qed\global\@qedtrue}%
\def\restore@qed{\global\let\if@qed\saved@ifqed}

\def\@Proof{%
   \par\removelastskip\bigskip\penalty100
   \save@set@qed
   \noindent\ProofFont{\TheWordProof\enskip}%
}%
\def\@Proofof#1{%
   \par\removelastskip\bigskip\penalty100
   \save@set@qed
   \noindent\ProofFont{\TheWordProofof{#1}\enskip}%
}%
\def\@endProof{%
   \qed\restore@qed
   \penalty-100 \medskip
}
\def\@endProofof{%
   \qed\restore@qed
   \penalty-100 \medskip
}
\def\@qed#1{%
\if@qed                                 
     \global\@qedfalse
        \ifmmode\ifinner\pushright{#1}
        \else\eqno{\qedsymbol}\fi
        \else\pushright{#1}\fi%
\else\ifhmode\ifinner\else\par\fi\fi
\fi}
\def\@pushright#1{%
  {\ifvmode                             
       \null\hfill{#1}\par              
  \else\ifmmode\maths@pushright{\hbox{#1}}
       \else\ifinner\hbox@pushright{#1}
            \else\parag@pushright{#1}
  \fi  \fi  \fi
}}%
\def\maths@pushright#1{{%
  \ifinner
     \hbox@pushright{#1}%
  \else
     \eqno#1
     \def\]{$$\ignorespaces}
  \fi
}}%
\usepackage{pictex}
\usepackage{amsmath,amsfonts,amscd}
\usepackage{amssymb}
\usepackage{theorem}
\usepackage{mathrsfs}
\usepackage{latexsym}
\usepackage{euscript}

{\theorembodyfont{\slshape}
\newtheorem{theorem}{Theorem}[section]
\newtheorem{proposition}[theorem]{Proposition}
\newtheorem{lemma}[theorem]{Lemma}
\newtheorem{corollary}[theorem]{Corollary}
}
{\theorembodyfont{\rmfamily}

\newtheorem{remark}[theorem]{Remark}

}
\newcommand\cuadrado{{\unskip\nobreak\hfil\penalty50\hskip2em\vadjust{}
\nobreak\hfil$\Box$\parfillskip=0pt\finalhyphendemerits=0\par}}

\def\N{\mathbb{N}}
\def\K{\mathbb{K}}
\def\Z{\mathbb{Z}}
\def\R{\mathbb{R}}

\def\C{\mathbb{C}}

\def\proj{\mathbb{P}}

\def\dist{{\sf dist}}

\renewcommand{\a}{\alpha}

\renewcommand{\d}{\delta}

\newcommand{\e}{\varepsilon}

\def\DD{{\mathcal D}}

\def\scC{{\mathscr C}}
\def\scD{{\mathrm{deg}}}
\def\mZ{{\mathcal Z}}

\newcommand{\x}{\times}
\newcommand{\<}{\langle}
\renewcommand{\>}{\rangle}

\renewcommand{\bar}{\overline}
\renewcommand{\tilde}{\widetilde}
\renewcommand{\bar}{\overline}

\def\Oh{{\cal O}}

\newcommand{\Sing}{\mathrm{Sing}}

\newcommand{\bfE}{\mathop{\mathbf E}}
\newcommand{\Avr}{\mathop{\mathbf {Avr}}}
\newcommand{\Prob}{\mathop{\rm Prob}}
\newcommand{\bfd}{\mathbf d}
\newcommand{\Hd}{{\mathcal H}_{\bfd}}

\newcommand{\cond}{{\sf cond}}

\def\algorithm{\begin{center}
               \begin{minipage}{6in}
               \begin{tabbing}
               \marks}
\def\falgorithm{\end{tabbing}
                \end{minipage}
                \end{center}}
\def\marks{nn\= nn\= nn\= nn\= nn\= nn\= nn\= \kill}

\def\PR{{\rm P}_{\kern-1pt\R}}
\def\PC{{\rm P}_{\kern-1pt\C}}
\def\NPR{{\rm NP}_{\kern-1pt\R}}
\def\NPC{{\rm NP}_{\kern-2pt\C}}
\def\coNPC{{\rm coNP}_{\kern-2pt\C}}
\def\coNPR{{\rm coNP}_{\kern-2pt\R}}
\def\DNPR{{\rm DNP}_{\kern-1pt\R}}
\def\DNPC{{\rm DNP}_{\kern-2pt\C}}
\def\PAR{{\rm PAR}_{\kern-1pt\R}}
\def\PARC{{\rm PAR}_{\kern-1pt\C}}
\def\PHR{{\rm PH}_{\kern-1pt\R}}
\def\PHC{{\rm PH}_{\kern-1pt\C}}
\def\DPHR{{\rm DPH}_{\kern-1pt\R}}
\def\DPHC{{\rm DPH}_{\kern-1pt\C}}

\def\FPR{{\rm FP}_{\kern-1pt\R}}
\def\FPC{{\rm FP}_{\kern-1pt\C}}
\def\FPAR{{\rm FPAR}_{\kern-0.4pt\R}}
\def\FPARC{{\rm FPAR}_{\kern-0.4pt\C}}

\def\CPRi{{\rm \#P}_{\kern-2pt\R}}
\def\CPRd{{\rm D\#P}_{\kern-2pt\R}}
\def\CPCi{{\rm \#P}_{\kern-2pt\C}}
\def\CPCd{{\rm D\#P_{\kern-2pt\C}}}
\def\gCPCi{{\rm \#P}^{\ast}_{\kern-2pt\C}}
\def\gCPRi{{\rm \#P}^{\ast}_{\kern-2pt\R}}

\def\FPK{{\rm FP}\kern-1.5pt_K}
\def\FPL{{\rm FP}\kern-1.5pt_L}
\def\CPK{\#{\rm P}\kern-1.5pt_K}
\def\CPL{\#{\rm P}\kern-1.5pt_L}
\def\CHNK{\#\mbox{\sc HN}\kern-1.5pt_K}
\def\CHNL{\#\mbox{\sc HN}\kern-1.5pt_L}
\def\PK{{\rm P}\kern-1.5pt_K}
\def\PL{{\rm P}\kern-1.5pt_L}
\def\NPK{{\rm NP}\kern-1.5pt_K}
\def\NPL{{\rm NP}\kern-1.5pt_L}

\def\FEASR{{\mbox{\sc Feas}_{\kern-0.5pt\R}}}
\def\FEASC{{\mbox{\sc Feas}_{\kern-1pt\C}}}
\def\FEASRbit{{\mbox{\sc Feas}^{\Z}_{\kern-1pt\R}}}
\def\FEASCbit{{\mbox{\sc Feas}^{\Z}_{\kern-1pt\C}}}
\def\HNC{{\mbox{\sc HN}_{\kern-1pt\C}}}
\def\PHNC{{\mbox{\sc PHN}_{\kern-1pt\C}}}
\def\HNCbit{{\mbox{\sc HN}^{\Z}_{\kern-1pt\C}}}
\def\QASC{{\mbox{\sc QAS}_{\kern-1pt\C}}}
\def\QASCbit{{\mbox{\sc QAS}^{\Z}_{\kern-1pt\C}}}
\def\PQAS{{\mbox{\sc ProjQAS}_{\kern-1pt\C}}}
\def\BPQAS{{\mbox{\sc BiProjQAS}_{\kern-1pt\C}}}

\def\DIMR{{\mbox{\sc Dim}_{\kern-0.5pt\R}}}
\def\DIMC{{\mbox{\sc Dim}_{\kern-0.5pt\C}}}
\def\DIMadd{{\mbox{\sc Dim}_{\kern-0.5pt\add}}}
\def\DIMRbit{{\mbox{\sc Dim}^{\Z}_{\kern-0.5pt\R}}}
\def\DIMCbit{{\mbox{\sc Dim}^{\Z}_{\kern-0.5pt\C}}}

\def\REACH{{\mbox{\sc Reach}_{\kern-0.5pt\R}}}
\def\REACHbit{{\mbox{\sc Reach}^{\Z}_{\kern-0.5pt\R}}}
\def\CREACHbit{{\mbox{\sc CReach}^{\Z}_{\kern-0.5pt\R}}}

\def\GapHNC{{\mbox{\sc $\Delta$HN}_{\kern-1pt\C}}}
\def\GapHNCbit{{\mbox{\sc $\Delta$HN}^{\Z}_{\kern-1pt\C}}}

\def\gGapC{{\mbox{\sc Gap}^\ast_{\kern-2pt\C}}}
\def\gGapR{{\mbox{\sc Gap}^\ast_{\kern-2pt\R}}}
\def\GapC{{\mbox{\sc Gap}_{\kern-1pt\C}}}

\newcommand{\HHNC}{\mbox{\sc HN}_{\kern-1pt\C}}

\def\Pp{{\proj^p}}
\def\Pm{{\proj^m}}
\def\Ppm{{\proj^{p-m}}}
\def\munorm{\mu_{\rm norm}}
\def\disc{{\sf disc}}

\def\Hi{{\mathrm{H}}}

\begin{document}
\begin{title}
{\Large {\bf Smoothed analysis of complex conic condition numbers}}
\end{title}
\author{Peter B\"urgisser\thanks{Dept.\ of Mathematics, University of
Paderborn, Germany. Partially supported
by DFG grant BU 1371.},
Felipe Cucker\thanks{Dept.\ of
Mathematics, City University of Hong Kong,
Kowloon Tong, Hong Kong. Partially supported by
CityU SRG grant 7001860.}, and Martin Lotz$^{\dagger}$}
\date{}
\makeatletter
\maketitle

\begin{quote}{\small
{\bf Abstract.} Smoothed analysis of complexity bounds and condition
numbers has been done, so far, on a case by case basis. In this
paper we consider a reasonably large class of condition numbers
for problems over the complex numbers and we obtain smoothed
analysis estimates for elements in this class depending only on
geometric invariants of the corresponding sets of ill-posed inputs.
These estimates are for a version of smoothed analysis
proposed in this paper which, to the best of our knowledge,
appears to be new.
Several applications
to linear and polynomial equation solving
show that estimates obtained in this way
are easy to derive and quite accurate.
}\end{quote}


\section{Introduction}\label{se:intro}

\subsection{Conic condition numbers---Main results}\label{sec:gf}

A distinctive feature of the computations considered in numerical
analysis is that they are affected by errors. A main character in
the understanding of the effects of these errors is the {\em
condition number} of the input at hand. This is a positive number
which, roughly speaking, quantifies the effects just mentioned when
computations are performed with infinite precision but the input has
been modified by a small perturbation. It depends only on the data
and the problem at hand. The best known condition number is that for
matrix inversion and linear equation solving. For a square matrix
$A$ it takes the form $\kappa(A)=\|A\|\|A^{-1}\|$
and was independently introduced by Goldstine and
von Neumann~\cite{vNGo47} and Turing~\cite{Turing48}.

Condition numbers occur in endless instances of round-off analysis.
They also appear as a parameter in complexity bounds for a variety
of iterative algorithms. Yet, condition numbers are not easily
computable. It has even been conjectured~\cite{Renegar94} that
computing the condition number $\scC(a)$ for a certain data $a$ is
at least as difficult as solving the problem for which $a$ is a
data.
A way out for this situation is to assume a probability measure on
the set of data and to study the condition number of this data as a
random variable.

The above ideas have been systematized in a number of places.
Notably, Blum~\cite{Blum90} suggested a complexity theory for
numerical algorithms parameterized by a condition number $\scC(a)$
of input data (in addition to input size). Then,
Smale~\cite[\S1]{Smale97} extended this suggestion by proposing to
obtain estimates on the probability distribution of $\scC(a)$.
Combining both ideas, he argued, one can give probabilistic bounds 
on the complexity of numerical algorithms.

Classically, probabilistic analysis of condition numbers takes two
forms: bounds on the tail of the distribution of $\scC(a)$
---showing that it is unlikely that $\scC(a)$ will be large--- and
bounds on the expected value of $\ln(\scC(a))$ ---estimating the
average loss of precision and average running time---. Examples of
such results abound for a variety of condition
numbers~\cite{ChC01,CW01,Demmel88,Edelman88,Bez2,TTY98}.

Recently D.~Spielman and S.-H.~Teng~\cite[\S3]{ST:02}
suggested a new approach
to Smale's agenda above. The idea (e.g., for the distribution's tail)
is to replace showing that
\begin{quote}
``it is unlikely that $\scC(a)$ will be large''
\end{quote}
by showing that
\begin{quote}
``for all $a$ and all slight random perturbation $\Delta a$, it
is unlikely that $\scC(a+\Delta a)$ will be large.''
\end{quote}
A survey of this approach, called {\em smoothed analysis},
can be found in~\cite{ST:02}. We briefly describe its main
features in \S\ref{sec:SA}.

The goal of this paper is to give bounds for the smoothed
analysis (both tail and expected value) for a large class
of condition numbers for problems over the complex numbers.
We assume our data space is $\C^{p+1}$, endowed with a
Hermitian product $\langle\ ,\ \rangle$. We say that $\scC$ is a
{\em conic condition number} if there exists an
algebraic cone $\Sigma\subset \C^{p+1}$ (the set of {\em
ill-posed inputs}) such that, for all data $a$,
\begin{equation*}
    \scC(a)=\frac{\|a\|}{\dist(a,\Sigma)},
\end{equation*}
where $\|\ \|$ and $\dist$ are the norm and distance
induced by $\langle\ ,\ \rangle$, respectively.

As defined above, $\kappa(A)$ is not conic since the operator
norm $\|\ \|$ is not induced by a Hermitian product. Replacing
this norm by the Frobenius norm $\|\ \|_F$ yields the
(commonly considered) version 
$\kappa_F(A) := \|A\|_F \|A^{-1}\|$ of $\kappa(A)$.
The Condition Number Theorem of Eckart and Young~\cite{EckYou}
then states that $\kappa(A)_F$ is conic
(with $\Sigma$ the set of singular matrices). 
Other examples
can be found in~\cite{ChC05}, where a certain property
(related with the so called level-2 condition numbers) is proved
for conic condition numbers. Furthermore, it is argued 
in~\cite{Demmel87} that for many problems, their condition number 
can be bounded by a conic one. 

Note that, since $\Sigma$ is a cone, for all $z\in\C\setminus\{0\}$,
$\scC(a)=\scC(za)$. Hence, we may restrict to
data $a\in\Pp:=\proj^p(\C)$ for which the condition number
takes the form
\begin{equation}\label{eq:conic-con}
  \scC(a)=\frac{1}{d_{\Pp}(a,\Sigma)}
\end{equation}
where, abusing notation, $\Sigma$ is interpreted now as a subset
of $\Pp$ and $d_{\Pp}$ denotes the projective distance in $\Pp$
(precise definitions follow in \S\ref{ss:prelim} below).
We will denote by $B(a,\sigma)$ the open ball of radius $\sigma$
around $a$ in $\Pp$ with respect to projective distance.

In what follows we assume $\Sigma$ is purely dimensional
and we write $m=\dim(\Sigma)$. Recall that the degree
$\scD(\Sigma)$ of $\Sigma$ equals (cf.~\cite{shaf:74})
\begin{equation*}
  \scD(\Sigma)=\min\{\ell\mid \#(\Sigma \cap\proj^{p-m})\leq \ell
      \mbox{ for almost all $\proj^{p-m}\subset\Pp$}\}.
\end{equation*}

Our main result is the following.

\begin{theorem}\label{th:smoothedCN}
Let $\scC$ be a conic condition number with set of ill-posed inputs
$\Sigma\subset\Pp$, of pure dimension $m$, $0< m< p$.
Then, for all $a\in\Pp$, all
$\sigma\in(0,1]$, and all $t\geq \frac{p\sqrt2}{p-m}$,
we have
\begin{equation*}
   \Prob_{z\in B(a,\sigma)}\{\scC(z)\geq t\}\leq
  K(p,m)\scD(\Sigma)
  \left(\frac{1}{t\sigma}\right)^{2(p-m)}
   \left(1+\frac{p}{p-m}\frac{1}{t\sigma}\right)^{2m}
\end{equation*}
and
\begin{equation*}
   \bfE_{z\in B(a,\sigma)}(\ln\scC(z))\leq
  \frac{1}{2(p-m)}\left(\ln K(p,m)+\ln\scD(\Sigma)
   +3\right)+\ln\frac{pm}{p-m}+2\ln\frac1\sigma,
\end{equation*}
with the constant $K(p,m):=2\frac{p^{3p}}{m^{3m}(p-m)^{3(p-m)}}$.
\end{theorem}

We will devote \S\ref{sec:appl} to derive applications
of Theorem~\ref{th:smoothedCN} to some condition numbers
which occur in the literature.

In most of our applications, the set of ill-posed inputs $\Sigma$
is a hypersurface. That is,
$\Sigma$ is the zero set $\mZ(f)$ of a nonzero homogeneous 
polynomial~$f$
and thus $\scD(\Sigma)$ is at most the degree of $f$.
In this case, we have the following easy to apply corollary.


\begin{corollary}\label{cor:hypersurface}
Let $\scC$ be a conic condition number with set of ill-posed inputs
$\Sigma\subseteq\Pp$. Assume $\Sigma\subseteq\mZ(f)$ with
$f\in\C[X_0,\ldots,X_p]$ homogeneous of degree $d$.
Then, for all $a\in\Pp$, all
$\sigma\in(0,1]$, and all $t\geq p\sqrt2$,
\begin{equation*}
   \Prob_{z\in B(a,\sigma)}\{\scC(z)\geq t\}\leq
    2p^3e^3 d
   \left(\frac{1}{t\sigma}\right)^{2}
   \left(1+p\frac{1}{t\sigma}\right)^{2(p-1)}
\end{equation*}
and
\begin{equation*}
   \bfE_{z\in B(a,\sigma)}(\ln\scC(z))\leq
  \frac{7}{2}\ln p+\frac12\ln d
   +4+2\ln\frac1\sigma.
\end{equation*}
\end{corollary}


The main idea towards the proof of Theorem~\ref{th:smoothedCN}
is to reformulate
the probability distribution of a conic condition number
as a geometric problem in a complex projective space.
Indeed, for $V\subseteq \Pp$ we denote
by $v(V)$ the volume of $V$, and by $V_\e$ the
$\e$-tube around $V$ in $\Pp$
(precise definitions follow in \S\ref{ss:prelim} below).
With this notation,
\begin{equation*}
  \Prob_{z\in B(a,\sigma)}\left\{\scC(z)
  \geq \frac{1}{\e}\right\}
  =\Prob_{z\in B(a,\sigma)}\left\{d_{\Pp}(z,\Sigma)\leq \e\right\}=
  \frac{v(\Sigma_\e\cap B(a,\sigma))}{v(B(a,\sigma))}.
\end{equation*}
The first claim in Theorem~\ref{th:smoothedCN} will thus follow
from the following purely geometric statement.

\begin{theorem}\label{th:tubos}
Let $V$ be a projective variety in $\Pp$ of
pure dimension $0<m<p$. Moreover, let $a\in\Pp$,
$\sigma\in(0,1]$, and $0<\e\leq \frac1{\sqrt{2}}\frac{p-m}{p}$.
Then we have
$$
  \frac{v(V_\e\cap B(a,\sigma))}{v(B(a,\sigma))}
  \leq K(p,m)\, \scD(V)
  \left(\frac{\e}{\sigma}\right)^{2(p-m)}
  \left(1+\frac{p}{p-m}\frac{\e}{\sigma}\right)^{2m}.
$$
\end{theorem}

One of the central tools in the derivation
of Theorem~\ref{th:tubos} is integral geometry. An essential
formula of integral geometry~\cite[\S15.2]{sant:76} allows
to relate the volume of certain geometric objects to the expected
volume of their intersection when they are moved at random.
A simple application is the equality
$v(V)=\scD(V) v(\proj^m)$ for the volume of
an irreducible $m$-dimensional subvariety $V\subseteq \proj^p$. In
order to obtain a corresponding bound for
$V_\e\cap B(a,\sigma)$, a more
sophisticated use of this equality is needed
(cf.~Lemma~\ref{lem:3}).


\subsection{Relation to previous work}\label{sec:SA}

Let $\K=\R$ or $\C$.
In the study of the behaviour of a function $f\colon\K^n\to\R_+$
(e.g., a condition number, a complexity bound) two frameworks
have been
extensively used: worst-case and average-case. Recently, a third
framework has been proposed which arguably blends the best of the
former two. The worst-case framework studies the quantity
\begin{equation}\label{eq:worst}
    \sup_{a\in\K^n} f(a)
\end{equation}
and the average-case the quantity
\begin{equation}\label{eq:avg}
  \bfE_{z\in \Psi} f(z) =\int_{z\in\K^n} f(z)\psi(z)dz
\end{equation}
where $z\in \Psi$ means that the expected value is taken for a
random $z$ whose distribution $\Psi$ has density function $\psi$.
The {\em smoothed analysis} of $f$ studies the behaviour of
\begin{equation}\label{eq:sa1}
   \sup_{a\in\K^n}\bfE_{z\in N^n(0,\sigma^2)} f(a+z)
\end{equation}
(possibly for sufficiently small $\sigma$) where $N^n(0,\sigma^2)$
denotes the $n$-dimensional Gaussian
distribution over $\K$ with mean 0 and variance $\sigma^2$.
Note that while (\ref{eq:worst}) and (\ref{eq:avg})
usually yield functions on $n$, (\ref{eq:sa1}) yields a
function on $n$ and $\sigma$. It has been argued that smoothed
analysis interpolates between worst and average cases since it
amounts to the first for $\sigma=0$ and it approaches the second for
large $\sigma$. Instances of smoothed analysis can be found
in~\cite{CDW:05,DST,ST:02,ST:03,ST:04,Wsch:04}.

When $f$ is homogeneous of degree 0  ---e.g., a conic condition
number--- it makes sense to restrict $f$ to the projective space
$\proj^{n-1}(\K)$. In this case, it also makes sense to replace the
distribution $a+N^n(0,\sigma^2)$ by the uniform distribution
supported on the disk
$B(a,\sigma)\subseteq\proj^{n-1}$
and consider, instead of (\ref{eq:sa1}), the
following quantity
\begin{equation}\label{eq:sa2}
   \sup_{a\in \proj^{n-1}}\bfE_{z\in B(a,\sigma)} f(z).
\end{equation}
Note that in this case, the interpolation mentioned above is
transparent. When $\sigma=0$ the expected value amounts to $f(a)$
and we obtain worst-case analysis, while if $\sigma=1$
(the diameter of $\proj^{n-1}$)
the expected value is independent of $a$ and we obtain
average-case analysis.

It is this version of smoothed analysis we deal with in
this paper. To the best of our knowledge it appears here
for the first time. Note that while, technically, this
``uniform smoothed analysis'' differs from the Gaussian
one considered so far, both share the viewpoint described
in~\S\ref{sec:gf} above.


We have already mentioned the
references~\cite{CDW:05,DST,ST:02,ST:03,ST:04,Wsch:04}
as instances of previous work in smoothed analysis.
In all these cases, an {\em ad hoc} argument is used to obtain
the desired bounds. This is in contrast with the goal
of this paper which is to provide general estimates
which can be applied to a large class of condition numbers.
We believe the applications in \S\ref{sec:appl} give
substance to this goal.

The idea of reformulating probability distributions as quotients
of volumes in projective spaces (or spheres) to estimate 
condition measures goes back at least 
to Smale~\cite{Smale81} and Renegar~\cite{Ren87a}. 
In particular,~\cite{Ren87a} uses this idea to show bounds 
on the probability distribution of a certain random variable in the 
average-case analysis of the complexity of Newton's method. 
Central to his argument is the 
fact that this random variable can be bounded by a conic 
condition number. The set of ill-posed inputs in~\cite{Ren87a} 
is a hypersurface. An extension of these results to the case 
of codimension greater than one was done by Demmel~\cite{Demmel88} 
where, in addition, an average-case analysis of several 
conic condition numbers is performed. Our paper is an 
extension of these arguments to the smoothed-analysis 
framework.

In a recent paper, Beltr\'an and Pardo~\cite{BePa:05} 
obtained estimates similar to those proved by Demmel 
(always for the average-case setting) when the input data 
$a$ is assumed to belong to a complex 
projective variety $V\subseteq\Pp$ and averages are taken
for the uniform distribution on $V$. 
An extension of Theorem~\ref{th:smoothedCN} in this direction  
is certainly doable,
but we have not included it in this paper.

Probably the most important extension of the present paper would 
be to obtain a result akin to Theorem~\ref{th:smoothedCN} 
(or Corollary~\ref{cor:hypersurface}) for
problems defined over the real numbers.  For the average-case 
setting Demmel~\cite{Demmel88}  
states such results. Unfortunately, his results directly
rely on an unpublished report by Ocneanu dating from 1985, which
apparently contains an upper bound on the volume of tubes around a
real variety in terms of degrees (cf.~Theorem 4.3 in~\cite{Demmel88}). 
We are currently working towards an extension to the real case.

\section{Proof of Theorem~\ref{th:smoothedCN}}\label{sec:proofs}

\subsection{Distances and volumes in projective space}
\label{ss:prelim}

We refer to~\cite[Chapter~12]{bcss:95}
for a more detailed introduction to the concepts needed here.
A general reference for complex analytic geometry 
is~\cite{grha:78}.

The complex projective space $\proj^p:=\proj^p(\C)$ is defined 
as the set of
one dimensional complex subspaces of $\C^{p+1}$.
The space $\proj^p$ carries the structure of a compact
$2p$-dimensional real manifold.
A Hermitian inner product $\langle\ ,\ \rangle$ on
$\C^{p+1}$ induces a Riemannian distance $d_R$ on $\proj^p$
(called Fubini-Study distance), which is defined as
\begin{equation*}
 d_R(x,y)=\arccos \frac{|\langle
\overline{x},\overline{y}\rangle|}{\|\bar{x}\|\, \|\bar{y}\|}
 \quad\mbox{ for $x,y\in\Pp$},
\end{equation*}
where $\bar{x},\bar{y}$ are representatives of $x$ and $y$ in
$\C^{p+1}$, respectively, and $\|\ \|$ denotes the norm induced
by $\<\ ,\ \>$.

The natural projection
$\R^{2p+2}\setminus\{0\}\cong
\C^{p+1}\setminus\{0\}\rightarrow \proj^p$
factors through a (everywhere regular) projection
$\pi\colon S^{2p+1}\rightarrow \proj^p$
with fiber $S^1$.
It is easy to check that the restriction of the derivative $d\pi(x)$
to the orthogonal complement
of its kernel is orthogonal with respect to the Riemannian metrics
on $S^{2p+1}$ and $\proj^p$ induced by $\langle\ ,\ \rangle$.
By means of the Co-Area formula \cite[p.~241]{bcss:95},
this observation allows to reduce the computation of integrals on
$\proj^p$ to the computation of integrals on $S^{2p+1}$. More
precisely, for any integrable function
$f\colon \proj^p\rightarrow\R$ and measurable
$U\subseteq \proj^p$ we have
\begin{equation}\label{eq:integrate}
  \int_U f d\proj^p=
  \frac{1}{2\pi}\int_{\pi^{-1}(U)}f\circ \pi\,dS^{2p+1},
\end{equation}
where $d\proj^p$ and $dS^{2p+1}$ denote the volume forms
induced by $\langle\ ,\ \rangle$.

For an open subset $U\subseteq M$ of an $m$-dimensional Riemannian
manifold, we write $v(U):=\int_U dM$ for the {\em $m$-dimensional volume}
of $U$, where $dM$ is the volume form on $M$ induced by the
Riemannian metric. In particular, using (\ref{eq:integrate}) we get
for the complex projective space
\begin{equation}\label{eq:volpn}
  v(\proj^p)=\frac{1}{2\pi} v(S^{2p+1})=\frac{\pi^p}{p!}.
\end{equation}

Instead of the Riemannian metric $d_R$ on $\proj^p$, we will be working
with the associated projective metric $d_{\Pp}$,
which is defined
as
\begin{equation*}
  d_{\Pp}(x,y)=\sin d_R(x,y).
\end{equation*}
Unless otherwise stated, this is the distance function we will be
using throughout this paper. The use of this distance function is
motivated by our applications. In fact, for a conic condition
number with ill-posed set $\Sigma\subseteq \C^{p+1}$, 
$d_{\Pp}(a,\Sigma)$ (recall our abuse of notation in the
introduction) just gives the normalized distance of a representative
of $a$ to $\Sigma$.

\begin{center}
   \input fig_smoothed_1.pictex
\end{center}

We denote by $B(x,\e)=B_\Pp(x,\e)$ the open ball of
radius $\e$ around $x$ in $\proj^p$ (with respect to $d_{\Pp}$),
and by $S^p(x,\e)$ the sphere of radius $\e$ around $x$.
For a subset $V\subseteq \proj^p$ we define the
{\em $\e$-tube} around $V$ in $\proj^p$ to be the open set
\begin{equation*}
  V_\e:=\{ x\in \proj^p \mid d_{\Pp}(x,V)<\e\}.
\end{equation*}
We will also use the notation $v_\e(V):=v(V_\e)$ for the
volume of an $\e$-tube in $\proj^p$ around a subset
$V\subseteq \proj^p$. If we wish to stress the ambient
space in which the tube is considered, we will write
$v_{\e}^{\Pp}(V)$ instead. We will similarly do so
if the ambient space is a sphere.

For a purely $m$-dimensional subvariety $V\subseteq \proj^p$, 
the set $V\backslash \Sing(V)$ (where $\Sing(V)$ denotes 
the singular locus of $V$)
is a real $2m$-dimensional Riemannian manifold (with the metric
induced from $\proj^p$), and we define the volume of $V$ as
$v(V):=v(V\backslash \Sing(V))$. This coincides with any other
reasonable notion of volume.

\begin{lemma}\label{lem:1}
Let $\Ppm \subseteq \Pp$ and let $0<\e\leq 1$. Then
\begin{equation*}
 v_{\e}^{\Pp}(\Ppm)\leq v(\Ppm)v(\Pm) \e^{2m},
\end{equation*}
with equality if and only if $p-m=0$. In particular,
for the volume of a ball of radius $\e$ around
$x\in \Pp$ we have
\begin{equation*}
  v\left(B_{\Pp}(x,\e)\right)= v(\Pp)\e^{2p}.
\end{equation*}
\end{lemma}

\begin{Proof}
A ball of radius $\e$ in $\Pp$ with respect to $d_{\Pp}$ corresponds
to a ball in $\Pp$ of radius $\d=\arccos(\e)$ with respect to $d_R$.
>From Equation (\ref{eq:integrate}) we get the
identity
\begin{equation}\label{eq:1}
  v_{\e}^{\Pp}\left(\Ppm\right)
  =\frac{v_\delta^{S^{2p+1}}\left(S^{2(p-m)+1}\right)}{2\pi}.
\end{equation}
Recall that on the sphere we use the usual Riemannian metric induced from
the ambient space. We have thus reduced our problem to that of computing the
volume of a tube around a subsphere of a sphere. Expressions for
this volume are straightforward to calculate: for a sphere
$S^m\subseteq S^p$ we have
\begin{equation*}
   v_{\delta}^{S^{p}}\left(S^{m}\right)
   =v(S^m)v(S^{p-m-1})\int_0^{\delta}\cos(t)^{m}\sin(t)^{p-m-1}dt.
\end{equation*}
Plugging this into Equation~(\ref{eq:1}) we get
\begin{align*}\label{eq:2}
  v_{\e}^{\Pp}\left(\Ppm\right)
  &=\frac{v\left(S^{2(p-m)+1}\right)v\left(S^{2m-1}\right)}
     {2\pi}\int_0^{\delta}\cos(t)^{2(p-m)+1}\sin(t)^{2m-1}dt\\
  &\stackrel{(\ref{eq:volpn})}{=}
  2\pi\, v(\Ppm)v(\proj^{m-1})\int_0^{\delta}\cos(t)^{2(p-m)+1}
     \sin(t)^{2m-1}dt\\
  &=2\pi\, v(\Ppm)v(\proj^{m-1})\int_0^{\e}(1-u^2)^{p-m}u^{2m-1}du,
\end{align*}
where in the last step we used the substitution $u=\sin(t)$. For
$0< u\leq 1$ we have $(1-u^2)^{p-m}\leq 1$, with equality if and
only if $p-m=0$. Substituting this bound in the above equation and
evaluating the integral, we get
\begin{equation*}
v_{\e}^{\Pp}\left(\Ppm\right)\leq
  \frac{2\pi\, v(\Ppm)v(\proj^{m-1})}{2m}\e^{2m}
  =v(\Ppm)v(\Pm)\e^{2m},
\end{equation*}
where we used the fact that $v(\Pm)=v(\proj^{m-1}) \pi/m$ for the
last equality.
\end{Proof}

\subsection{A fact from integral geometry}\label{ss:int-geo}

We will repeatedly use a variation of a classical
formula from integral geometry.
Let $M,N\subseteq \proj^p$ be submanifolds of (real) dimension 
$2m$ and $2n$, respectively. The unitary group $G:=U(p+1)$
acts transitively on $\proj^{p}$ in a straightforward way.
A key result in integral geometry states that the expected 
volume of the intersection of $M$ with a
random translate $gN$ of $N$ satisfies
\begin{equation}\label{eq:prior}
  \frac{\bfE_{g\in G}(v(M\cap
  gN))}{v(\proj^{m+n-p})}=\frac{v(M)v(N)}{v(\proj^m)v(\proj^n)}.
\end{equation}
Hereby the expectation is taken with respect to the normalized
Haar measure on~$G$.
The above equality also holds if $M$ and $N$ are (possibly
singular) subvarieties of $\Pp$.
Equation~(\ref{eq:prior}) is easily derived, 
using~(\ref{eq:integrate}),
from the corresponding statement in~\cite[\S15.2]{sant:76} 
for spheres.

\subsection{Estimating the volume of patches of projective varieties}
\label{ss:vol-tubes}

The following lemma allows to estimate the volume of the intersection of a
projective variety~$V$ with a ball in terms of the degree of~$V$ and the radius of the ball.

\begin{lemma}\label{lem:3}
Let $V\subset\Pp$ be an irreducible $m$-dimensional projective
variety, $a\in\Pp$, $0<\e\leq 1$ and $V'=V\cap B_{\Pp}(a,\e)$.
Then
\begin{equation*}
  \frac{v(V')}{v(\Pm)}\leq
  \scD(V) {p\choose m}\e^{2m}.
\end{equation*}
\end{lemma}

\begin{Proof}
Taking $M=\proj^{p-m}$ and $N=V'$ in~(\ref{eq:prior})
we obtain
\begin{equation*}
  \frac{v(V')}{v(\Pm)}= \bfE(|gV'\cap\Ppm|)
\end{equation*}
where the expectation is over all $g$ in the unitary group $U_{p+1}$
taken w.r.t.\ the normalized Haar measure (so that $U_{p+1}$ has
volume 1). Since $|gV'\cap\Ppm|\leq |gV\cap\Ppm| \leq \scD(V)$ for
almost all $g\in U_{p+1}$ we obtain
\begin{equation*}
  \bfE(|gV'\cap\Ppm|)\leq \scD(V)
  \Prob_{g\in U_{p+1}}\{gV'\cap \Ppm\neq\emptyset\}.
\end{equation*}
Since $V'\subseteq B(a,\e)$ we have
\begin{equation*}
  \Prob_{g\in U_{p+1}}\{gV'\cap \Ppm\neq\emptyset\}
  \leq \Prob_{g\in U_{p+1}}\{gB(a,\e)\cap
  \Ppm\neq\emptyset\}
  =\frac{v_{\e}^{\Pp}(\Ppm)}{v(\Pp)}.
\end{equation*}
The statement now follows from Lemma~\ref{lem:1} using that
$v(\Pp)=\frac{\pi^p}{p!}$.
\end{Proof}

The following crucial lemma is the only step in our chain of argumentation
that fails to be true over $\R$.

\begin{lemma}{\bf \cite[Theorem~22]{BePa:05}}\label{lem:2}
Let $V\subset\Pp$ be an irreducible projective variety of
dimension $m\geq 1$, $y\in V$ and $0<\e\leq 1/\sqrt{2}$.
Then we have
\begin{equation}\tag*{\cuadrado}
  v(V\cap B_{\Pp}(y,\e))\geq \frac12v(\Pm)\e^{2m}.
\end{equation}
\end{lemma}
\medskip

\subsection{Bounding the expectation}\label{ss:expect}

The next result gives a convenient way to bound the expectation
of a nonnegative random variable whose tail probabilities can be
estimated by some power law.

\begin{proposition}\label{prop:technical}
Let $X$ be a nonnegative, absolutely continuous, random variable
and $\alpha,t_0,K$ be positive constants
satisfying $\Prob\{X\geq t\}\leq K t^{-\alpha}$ for
all $t\geq t_0$. Then we have
\begin{equation*}
 \bfE(\ln X)\leq \ln t_0+\frac1{\a}\left(\ln K + 1\right).
\end{equation*}
Moreover, if $t_0\leq K^{\frac{1}{\alpha}}$ then
$\bfE(\ln X)\leq \frac1{\a}\left(\ln K + 1\right)$.
\end{proposition}

\begin{Proof}
Define the monotonically decreasing function $g:(0,1)\to\R$ by
\begin{equation*}
   g(y)=\left\{\begin{array}{ll}
         -\frac1{\alpha}\ln\left(\frac y{K}\right) &
          \makebox{if $y\leq Kt_0^{-\alpha}$}\\
         \ln t_0 & \makebox{otherwise.}\end{array}\right.
\end{equation*}
We claim that $\Prob\{\ln X\geq g(y)\}\leq y$ for all $y\in(0,1)$.
Indeed, if $y\leq Kt_0^{-\alpha}$ then there exists $t\geq t_0$
such that $y=Kt^{-\alpha}$. Therefore,
\begin{equation*}
  g(y)=-\frac1{\alpha}\ln\left(\frac y{K}\right)
      = \ln t
\end{equation*}
and
\begin{equation*}
  \Prob\{\ln X\geq g(y)\}=\Prob\{\ln X\geq \ln t\}
  =\Prob\{X\geq t\}\leq Kt^{-\alpha}=y.
\end{equation*}
If, instead, $y> Kt_0^{-\alpha}$ then
\begin{equation*}
  \Prob\{\ln X\geq g(y)\}=\Prob\{\ln X\geq \ln t_0\}
  =\Prob\{X\geq t_0\}\leq Kt_0^{-\alpha}<y.
\end{equation*}
Using~\cite[Prop.~2, Ch.~11]{bcss:95} it follows that
\begin{eqnarray*}
  \bfE(\ln X)&\leq& \int_0^1 g(y) dy\\
  &=&-\int_0^{Kt_0^{-\alpha}}\frac1{\a}\ln(y/K)dy
     +\int_{Kt_0^{-\alpha}}^1 \ln t_0\, dy\\
  &\leq& -\int_0^1\frac1{\a}\ln(y/K)\,dy
     +\int_{0}^1 \ln t_0\, dy\\
  &=&\frac1{\a}y(\ln y-1)\biggl|_{1}^0
      +\frac1{\a}\ln K +\ln t_0\\
  &=&\frac1{\a}(1+\ln K) +\ln t_0.
\end{eqnarray*}
If $t_0\leq K^{\frac{1}{\alpha}}$ then $Kt_0^{-\alpha}\geq 1$
and the integral above has only its first term.
\end{Proof}

\subsection{Proof of main results}\label{ss:proof-main}

\proofof{Theorem~\ref{th:tubos}}
It is sufficient to prove the assertion for an irreducible $V$.
In order to see this recall that 
$\scD(V)=\scD(V_1)+\cdots+\scD(V_q)$,
where $V_1,\ldots,V_q$ are the irreducible components of $V$
which we assume to be all of the same dimension.

So we assume that $V$ is irreducible.
We follow the arguments in~\cite[Proof of Theorem~16]{BePa:05}.
Fix $\e_1\in(0,1]$ such that $0<\e_1-\e\leq\frac1{\sqrt2}$
(we will specify $\e_1$ later). For each $z\in V_\e$ there exists
$y\in V$ such that $d_{\Pp}(z,y)\leq \e$ and hence
$B(y,\e_1-\e)\subseteq B(z,\e_1)$.
\smallskip

\begin{center}
   \input fig_smoothed_2.pictex
\end{center}

Since $\e_1-\e\leq\frac1{\sqrt2}$ we may use Lemma~\ref{lem:2} to
obtain
\begin{equation}\label{eq:lower}
  v(V\cap B(z,\e_1))\geq
  v(V\cap B(y,\e_1-\e)) \geq
  \frac12 v(\Pm)(\e_1-\e)^{2m}.
\end{equation}
In order to estimate $v(V_\e\cap B(a,\sigma))$ we put
$V':=V\cap B(a,\sigma+\e_1)$ and note that
$V\cap B(z,\e_1)\subseteq V'$ for all $z\in V_\e\cap B(a,\sigma)$.
\smallskip

\begin{center}
   \input fig_smoothed_3.pictex
\end{center}

Using (\ref{eq:lower}) we have
\begin{eqnarray*}
 \frac{v(V_\e\cap B(a,\sigma))}{v(\Pp)}
 &=&
 \frac{1}{v(\Pp)}\int_{z\in V_\e\cap B(a,\sigma)}1\;dz \\
 &\leq&
 \frac{1}{v(\Pp)}\int_{z\in V_\e\cap B(a,\sigma)}
 \frac{2\; v(V\cap B(z,\e_1))}{v(\Pm)(\e_1-\e)^{2m}}dz \\
 &\leq&
 \frac{2}{v(\Pm)(\e_1-\e)^{2m}}\frac1{v(\Pp)}
 \int_{z\in \Pp} v(V'\cap B(z,\e_1))dz.
\end{eqnarray*}
In addition,
\begin{eqnarray*}
 \frac1{v(\Pp)}\int_{z\in \Pp} v(V'\cap B(z,\e_1))dz
 &=&\int_{g\in U_{p+1}} v(V'\cap B(gz_0,\e_1))dg\\
 &\stackrel{(\ref{eq:prior})}{=}&
  v(\Pm)\frac{v(V')}{v(\Pm)}\frac{v(B(z_0,\e_1))}{v(\Pp)},
\end{eqnarray*}
where $z_0$ is any point in $\Pp$ and the second equality
follows from~(\ref{eq:prior}).
Using Lemma~\ref{lem:1} we conclude that
$$
  \frac1{v(\Pp)}\int_{z\in \Pp} v(V'\cap B(z,\e_1))dz
  =v(V')\e_1^{2p}.
$$
On the other hand, by Lemma~\ref{lem:3}, we have
$$
  \frac{v(V')}{v(\Pm)}\leq \scD(V) {p\choose m}
  (\sigma+\e_1)^{2m}
$$
since $V'= V \cap B(a,\sigma+\e_1)$. Combining all the above
we get the estimate
$$
  \frac{v(V_\e\cap B(a,\sigma))}{v(\Pp)} \leq
  \frac{2}{(\e_1-\e)^{2m}}\frac{v(V')}{v(\Pm)}\e_1^{2p}
   \le \frac{2\e_1^{2p}}{(\e_1-\e)^{2m}} \scD(V) {p\choose m}
  (\sigma+\e_1)^{2m}.
$$
Using again $v(B(a,\sigma))=v(\Pp)\sigma^{2p}$ it follows that
$$
  \frac{v(V_\e\cap B(a,\sigma))}{v(B(a,\sigma))} \leq
  \frac{2}{(\e_1-\e)^{2m}} \left(\frac{\e_1}{\sigma}\right)^{2p}
  \scD(V) {p\choose m}(\sigma+\e_1)^{2m}.
$$
We finally choose $\e_1:=\frac{p}{p-m}\e$. Note that
then
$$
  \e_1-\e=\frac{m}{p-m}\e\leq \frac{1}{\sqrt2}
$$
as we needed, the inequality since
$\e\leq\frac1{\sqrt2}\frac{p-m}{m}$. We obtain
$$
  \frac{v(V_\e\cap B(a,\sigma))}{v(B(a,\sigma))} \leq
  \frac{2p^{2p}}{m^{2m}(p-m)^{2(p-m)}}{p\choose m}\scD(V)
  \left(\frac{\e}{\sigma}\right)^{2(p-m)}
   \left(1+\frac{p}{p-m}\frac{\e}{\sigma}\right)^{2m} .
$$
Taking into account the estimate
${p\choose m} \le \frac{p^p}{m^m(p-m)^{p-m}}$
(which readily follows from the binomial expansion of
$p^p = (m+ (p-m))^p$) we finish the proof.
\endproofof

\proofof{Theorem~\ref{th:smoothedCN}}
The inequality for the tail follows directly from
Theorem~\ref{th:tubos}. For the expectation estimate,
let $\e_0:=\frac{p-m}{pm}\sigma$ and $t_0:=\e_0^{-1}$.
Note that, for $\e\leq \e_0$,
$$
  \left(1+\frac{p}{p-m}\frac{\e}{\sigma}\right)^{2m}
  \leq \left(1+\frac{1}{m}\right)^{2m}\leq e^2
$$
and thus
$$
  \frac{v(V_\e\cap B(a,\sigma))}{v(B(a,\sigma))} \leq
  K(p,m)\scD(\Sigma)
  \left(\frac{\e}{\sigma}\right)^{2(p-m)} e^2.
$$
Therefore, for all $t\geq t_0$, writing $\e=1/t$,
\begin{eqnarray*}
  \Prob_{z\in B(a,\sigma)}\{\scC(z)\geq t\}
  &=&\Prob_{z\in B(a,\sigma)}\{d(z,\Sigma)\leq \e\}\\
  &=&\frac{v(V_\e\cap B(a,\sigma))}{v(B(a,\sigma))}\\
  &\leq&
  K(p,m)\scD(\Sigma)
  \left(\frac{1}{\sigma}\right)^{2(p-m)} e^2 t^{-2(p-m)}.
\end{eqnarray*}
A straightforward application of
Proposition~\ref{prop:technical} yields
$$
  \bfE_{z\in B(a,\sigma)}(\ln \scC(z)) \leq
  \frac{1}{2(p-m)}\left(\ln K(p,m)+\ln\scD(\Sigma)
   +3\right)+\ln\frac{pm}{p-m}+2\ln\frac1\sigma.
$$
\endproofof

\proofof{Corollary~\ref{cor:hypersurface}}
Put $\Sigma'=\mZ(f)$ and note that
$\scC(a)=\frac{1}{d_{\Pp}(a,\Sigma)} \le 
\frac{1}{d_{\Pp}(a,\Sigma')}$.
The assertion follows from Theorem~\ref{th:smoothedCN} 
applied to $\Sigma'$ and the inequality
\begin{equation}\tag*{\cuadrado}
K(p,p-1)=2\frac{p^{3p}}{(p-1)^{3p-3}}p
   = 2\left[\left(1+\frac{1}{p-1}\right)^{p-1}\right]^3
      p^3
  \leq 2e^3 p^3.
\end{equation}

\section{Some Applications}\label{sec:appl}

In this section we obtain smooth analysis estimates
for the condition numbers of four problems: linear equation
solving, Moore-Penrose inversion, eigenvalue computations, 
and polynomial equation solving. For the first two, instances 
of such analysis already exist and we therefore compare
our results with those in the literature. The following
differences, however, should be noted.
Firstly, these analyses were done for
problems over the reals. Secondly, they hold
within the Gaussian framework
for smoothed analysis described in~\S\ref{sec:SA}. The
first feature is not important since a cursory
look at the refered proofs shows that similar results hold
for complex matrices. One should though keep in mind the second.

\subsection{Linear equation solving}\label{ss:les}

The first natural application of our result is for
the classical condition number $\kappa(A)$.
In~\cite{Wsch:04}, M.~Wschebor showed (solving
a conjecture posed in~\cite{ST:02}) that, for
all $n\times n$ real matrices $M$ with $\|M\|\leq 1$,
all $0<\sigma\leq 1$ and all $t>0$
$$
   \Prob_{E\in N^{n^2}(0,\sigma^2)}(\kappa(M+E)\geq t)\leq
   \frac{Kn}{\sigma t}
$$
with $K$ a universal constant. Note that,
by Proposition~\ref{prop:technical}, this implies
$$
   \bfE_{E\in N^{n^2}(0,\sigma^2)}(\ln \kappa(M+E))\leq
   \ln n+\ln\frac1{\sigma}+\ln K +1.
$$
We next compare Wschebor's result with
what can be obtained from Corollary~\ref{cor:hypersurface}.
To do so, we first note that, for $A\in\C^{n\times n}$,
$$
  \kappa(A)=\|A\|\|A^{-1}\|\leq \|A\|_F\|A^{-1}\|=:\kappa_F(A)
$$
and that,
by the Condition Number Theorem of Eckart and
Young~\cite{EckYou}
(see also~\cite[Theorem~1, Chapter~11]{bcss:95}),
$\|A^{-1}\|=d_F(A,\Sigma)^{-1}$. Here $\|\ \|_F$ and $d_F$
are the Frobenius norm and distance in $\C^{n\times n}$
which are induced by the Hermitian product
$(A,B)\mapsto{\mathsf{trace}}(A B^*)$. It follows that
$\kappa_F(A)$ is conic. We can thus give upper bounds for
$\kappa_F(A)$ and they will hold as well for $\kappa(A)$.

\begin{proposition}\label{prop:cuad}
For all $n\geq 1$, $0<\sigma\le 1$,
and $M\in\C^{n\times n}$ we have
$$
   \bfE_{A\in B(M,\sigma)}(\ln \kappa_F(A))\leq
    \frac{15}{2}\ln n+2\ln\left(\frac1\sigma\right)+4,
$$
where the expectation is over all $A$ uniformly distributed
in the disk of radius $\sigma$ centered at $M$ in projective space
$\proj^{n^2-1}$ (recall that we always use the projective and not
the Riemannian distance).
\end{proposition}

\begin{Proof}
The variety $\Sigma$ of singular matrices is a hypersurface in
$\proj^{n^2-1}$ of degree $n$. We now apply
Corollary~\ref{cor:hypersurface}.
\end{Proof}

Note, the bound in Proposition~\ref{prop:cuad} is of the same
order of magnitude than Wschebor's, worse by just a constant
factor. On the other hand, its derivation from
Corollary~\ref{cor:hypersurface} is rather immediate.
We next extend this bound to rectangular matrices.

\subsection{Moore-Penrose inversion}\label{ss:MPi}

Let $\ell\geq n$ and consider the space $\C^{\ell\times n}$ of
$\ell\times n$ rectangular matrices. Denote by
$\Sigma\subset\C^{\ell\x n}$ the subset of rank-deficient matrices.
Let $A\not\in\Sigma$ and let $A^\dagger$ denote its Moore-Penrose
inverse (see, e.g.,~\cite{BG,CM}). The condition number of $A$ (for
the computation of $A^\dagger$) is defined as
$$
  \cond^{\dagger}(A)=\lim_{\e\rightarrow 0}
  \sup_{\|\Delta A\|_2\leq \e}
  \frac{\|(A+\Delta A)^{\dagger} - A^{\dagger}\|_2 \|A\|_2}
       {\|A^{\dagger}\|_2 \|\Delta A\|_2}.
$$
This is not a conic condition number but it happens to be close to
one. One defines $\kappa^{\dagger}(A)=\|A\|_2 \|A^\dagger\|_2$
and, since $\|A^\dagger\|_2=\dist(A,\Sigma)^{-1}$~\cite{GoLoan}, we
obtain
$$
  \kappa^{\dagger}(A)=\frac{\|A\|_2}{\dist(A,\Sigma)}.
$$
In addition (see~\cite[\S III.3]{Stewart}),
$$
  \kappa^{\dagger}(A)\leq\cond^{\dagger}(A)\leq
  \frac{1+\sqrt{5}}{2}\kappa^{\dagger}(A).
$$
Thus, $\ln(\cond^{\dagger}(A))$ differs from
$\ln(\kappa^{\dagger}(A))$
just by a small additive constant. As for square matrices,
$\kappa^{\dagger}(A)$ is not conic since the operator
norm is not induced by a Hermitian product in
$\C^{\ell\times n}$. But, again, we can bound
$\kappa^{\dagger}(A)$ by the conic
condition number
$\kappa_F^{\dagger}(A):=\|A\|_F\|A^{\dagger}\|$.

A smoothed analysis for $\kappa^{\dagger}(A)$ was performed
in~\cite{CDW:05}.
Computer experiments reported in that paper,
however, suggest that the exhibited bounds, while sharp
when $\ell$ is close to $n$, are not so for more elongated
matrices. Actually, an empirical average
$\Avr(\ln \kappa^{\dagger}(A))$
was computed for several pairs $(n,\ell)$ and matrices
of the form $A=M+\Delta$ with $M$ a fixed ill-posed matrix and
$\Delta$ a small perturbation. It was then
mentioned~\cite[\S7]{CDW:05} that ``one sees that when one
fixes $n$ and lets $\ell$ increase the quantity
$\Avr(\ln\kappa^{\dagger}(A))$ decreases. This is in contrast
with the behaviour of [our bound]. It appears that our methods
are not sharp enough to capture the behaviour of
$\bfE(\ln\kappa^{\dagger}(A))$.'' As we next see,
the bounds following from Theorem~\ref{th:smoothedCN}
capture this behaviour much better.

The bound shown in~\cite{CDW:05} is of the form
\begin{equation}\label{eq:boundCDW}
  \displaystyle\sup_{A\in\R^{\ell\times n}}
  \bfE_{E\in N^{\ell n}(0,\sigma)}(\ln\kappa^{\dagger}(A+E))
  \leq \Oh\left(\ln \ell +\ln\frac1\sigma\right).
\end{equation}
It depends on $\ell$ and
tends to $\infty$ when $\ell$ does so.
Our next result shows that for large~$\ell$,
the expected value above (now with respect to
uniform perturbations)
is bounded by an expression depending only on $n$ and $\sigma$.

\begin{proposition}\label{prop:MP}
For all $n\geq 1$ and $0<\sigma\le 1$  we have
$$
  \limsup_{\ell\to\infty}\sup_{M\in\proj^{\ell n-1}}
   \bfE_{A\in B(M,\sigma)}
   (\ln\kappa_{F}^{\dagger}(A))\leq
    \left(n+\frac32\right)\ln (n) +n\ln 2+ 2+
      (n+1)\ln\frac 1\sigma.
$$
\end{proposition}

\begin{Proof}
It is well known that (the image in $\proj^{n\ell-1}(\C)$ of)
$\Sigma$ is a projective variety of codimension $\ell-n+1$ and
degree ${\ell \choose n-1}$
(see~\cite[Examples~12.1 and 19.10]{Harris92}).
By Theorem~\ref{th:smoothedCN}, for all $M\in\Pp$
and $t\ge t_0=1$ 
\begin{eqnarray*}
   \Prob_{A\in B(M,\sigma)}\{\kappa_F^{\dagger}(A)\geq t\}
  &\leq& K(p,m)\scD(\Sigma)
  \left(\frac{1}{t\sigma}\right)^{2(p-m)}
   \left(1+\frac{p}{p-m}\frac{1}{t\sigma}\right)^{2m}\\
  &\leq& K(p,m)\scD(\Sigma)
  \left(\frac{1}{t\sigma}\right)^{2(p-m)}
   \left(\frac{2p}{\sigma(p-m)}\right)^{2m}
\end{eqnarray*}
with
$$
   p=\ell n-1, \qquad
   m=\ell n-\ell+n-2,
  \qquad\mbox{and}\qquad
  \scD(\Sigma)= {\ell\choose n-1}.
$$
Therefore, by Proposition~\ref{prop:technical},
\begin{eqnarray*}
   \bfE_{A\in B(M,\sigma)}(\ln\kappa_F^{\dagger}(A))&\leq&
   \frac{1}{2(p-m)}\ln \bigg(K(p,m)\scD(\Sigma)
   \left(\frac{2p}{\sigma(p-m)}\right)^{2m}+1\bigg)
   + \ln\frac{1}{\sigma}.
\end{eqnarray*}
We next bound the logarithms of the expressions inside the
parenthesis.

To bound the binomial coefficients we use the
following estimates (see~\cite[(1.4.5)]{Lint99})
$$
 \ln {p\choose m} \le \ln\frac{p^{p}}{m^{m}(p-m)^{(p-m)}}
 \le p\, \Hi\Big(\frac{m}{p}\Big),
$$
where $\Hi$ denotes the binomial entropy function defined by
$\Hi(z)=-z\ln z-(1-z)\ln(1-z)$ for $z\in (0,1)$.
Note that $\Hi$ is monotonically increasing on $(0,\frac12)$
and $\Hi(z)=\Hi(1-z)$ for $z\in (0,1)$.

It will be convenient to use the asymptotic notations
$f(n,\ell)\sim g(n,\ell)$ and $f(n,\ell)\lesssim g(n,\ell)$
to express that
$\lim\limits_{\ell\to\infty}\frac{f(n,\ell)}{g(n,\ell)} =1$ and
$\limsup\limits_{\ell\to\infty}\frac{f(n,\ell)}{g(n,\ell)} \le 1$, 
respectively.

We obtain
$$
 p\, \Hi\Big(\frac{m}{p}\Big) =  p\,\Hi\Big(\frac{p-m}{p}\Big)
 \le p\,\Hi\Big(\frac{\ell -n+1}{\ell n -1}\Big)
 \sim \ell n\,\Hi\Big(\frac1n\big) \le \ell (1+\ln n)
$$
using
$$
 \Hi\Big(\frac1n\Big) =
 \frac{1}{n}\ln n +\frac{n-1}{n}\ln \frac{n}{n-1} \le \frac1n (1 +\ln n).
$$
Hence $\ln K(p,m) \simeq 3\ell(1+\ln n)$.
Similarly,
$$
 \ln\scD(\Sigma) = \ln {\ell\choose n-1}
 \le \ell\, H\Big(\frac{n}{\ell}\Big)\lesssim n\ln\frac{\ell}{n}.
$$
Finally,
$$
 2m \ln \bigg(\frac{2p}{\sigma(p-m)}\bigg) \simeq
  2\ell n (\ln n + \ln\frac2\sigma).
$$
Therefore,
\begin{align*}
 \sup_{M\in\proj^{\ell n -1}} \bfE_{A\in B(M,\sigma)}&(\ln\kappa_F^{\dagger}(A))\\
 &\lesssim\;\frac{1}{2\ell}\left(
  3\ell (1+\ln n) +n\ln\left(\frac{\ell}{n}\right) + 2\ell n 
   \left(\ln n + \ln\frac{2}{\sigma}\right) \right)
  + \ln\frac{1}{\sigma} \\
 &\le\; \left(n+\frac32\right)\ln n 
      +(n+1)\ln\frac{1}{\sigma} + n\ln2 + 2,
\end{align*}
which shows the claim.
\end{Proof}

\begin{remark}
The bound in Proposition~\ref{prop:MP} is independent
of $\ell$. Yet, its dependance on $n$ is linear and
the term on $\ln\frac1\sigma$ is multiplied by a
factor $n$. This is too large a bound. We now note that
bounds such as (\ref{eq:boundCDW}) also follow from our
results. For a very short derivation, note that
if a matrix $A$ is rank deficient then, $\det(\overline{A})=0$
where $\overline{A}$ is the $n\times n$ matrix obtained by removing
all rows of $A$ with index greater than $n$. Therefore
$\Sigma\subseteq\overline{\Sigma}=\{A\in\C^{\ell\times n}\mid
\det(\bar{A})=0\}$. This implies that, if $\|A\|_F=1$
$$
  \kappa_F^{\dagger}(A)\leq
  \frac{1}{d_{\proj^{\ell n-1}}(A,\bar{\Sigma})}.
$$
Since $\bar{\Sigma}$ is a hypersurface of degree~$n$,
an immediate application of
Corollary~\ref{cor:hypersurface} yields
\begin{equation*}
   \sup_{M\in\proj^{\ell n-1}}
   \bfE_{A\in B(M,\sigma)}(\ln\kappa^{\dagger}(A))\leq
   \frac{7}{2}\ln \ell + 4\ln n
   +4+2\ln\left(\frac1\sigma\right).
\end{equation*}
For small $\ell$ (say, polynomially bounded in $n$) this
last bound is better than that in Proposition~\ref{prop:MP}.
We conjecture that an asymptotic bound of the form
$\Oh(\ln(n)+\ln(1/\sigma))$ holds.
\end{remark}

\subsection{Eigenvalue computations}\label{ss:eigen}

Let $A\in\C^{n\times n}$ and $\lambda\in\C$ be a simple
eigenvalue of $\C$. For any sufficiently small perturbation
$\Delta A$ there exists a unique eigenvalue $\tilde\lambda$
of $A+\Delta A$ close to $\lambda$. It is known~\cite{Kato:66}
that 
\begin{equation}\label{eq:eig}
  |\lambda-\tilde\lambda|\leq \|P\|\|\Delta A\| +\Oh(\|\Delta A\|^2)
\end{equation}
where $P\in\C^{n\times n}$ is the {\em projection matrix} given by
$$
   P=(y^Hx)^{-1}xy^H. 
$$
Here $x$ and $y$ are right and left eigenvectors associated
to $\lambda$, respectively, (i.e., satisfying 
$Ax=\lambda x$ and $y^HA=\lambda y^H$)
and $y^H$ is the transpose conjugate of $y$. Note that 
$y^Hx$ is a scalar. Furthermore, 
inequality~(\ref{eq:eig}) is sharp in the sense that the factor 
$\|P\|$ can not be decreased. We can
then define
$$
  \kappa(A,\lambda):=\left\{
   \begin{array}{ll}
   \|P\| & \mbox{if $\lambda$ is simple}\\
   \infty & \mbox{otherwise} \end{array}\right.
$$
and the (absolute) condition number of $A$ for eigenvalue
computations
$$
  \kappa_{\sf eigen}(A):=\max_{\lambda}
                         \kappa(A,\lambda),
$$
where the maximum is over all the eigenvalues $\lambda$ of $A$. 
Note that $\kappa_{\sf eigen}(A)$ is homogeneous of degree 0
in $A$. Also, the set $\Sigma$ where $\kappa_{\sf eigen}$ is
infinite is the set of matrices having multiple eigenvalues.
Finally, Wilkinson~\cite{Wilkinson72} proved that
\begin{equation}\label{eq:wilki}
   \kappa_{\sf eigen}(A)\leq
   \frac{\sqrt{2}\|A\|_F}{\dist(A,\Sigma)}.
\end{equation}
In~\cite{Demmel88}, Demmel used the fact that the right-hand
side of (\ref{eq:wilki}) is conic to obtain bounds on the
tail of $\kappa_{\sf eigen}(A)$ for random $A$. We next
use it to obtain smoothed analysis estimates.

\begin{proposition}\label{prop:eigen}
For all $n\geq 1$ and $M\in\C^{n\times n}$,
$$
   \bfE_{A\in B(M,\sigma)}(\ln \kappa_{\sf eigen}(A))\leq
    8\ln n+2\ln\frac1\sigma + 5.
$$
\end{proposition}

\proof
Let $\chi_A$ be the characteristic polynomial of $A$.
This is a monic polynomial of degree $n$ whose coefficient
of degree $i$ is a homogeneous polynomial of degree $n-i$
in the entries of $A$.
Clearly, $A$ has multiple eigenvalues if and only if $\chi_A$
has multiple roots. This happens if and only if the 
discriminant $\disc(\chi_A)$ of $A$ is zero. 
The discriminant $\disc(\chi_A)$ is a polynomial in the entries
of~$A$, which can be expressed in terms of the eigenvalues
$\lambda_1,\ldots,\lambda_n$ of $A$ as follows 
$$ 
  \disc(\chi_A) =
  \prod_{i<j} (\lambda_i-\lambda_j)^2. 
$$ 
Note that
$\alpha\lambda_1,\ldots,\alpha\lambda_n$ are the eigenvalues of
$\alpha A$, for $\alpha\in\C$. Hence 
$$ 
   \disc(\chi_{\alpha A}) =\prod_{i<j} 
  (\alpha\lambda_i - \alpha\lambda_j)^2 
  = \alpha^{n^2 - n}\prod_{i<j} (\lambda_i-\lambda_j)^2. 
$$ 
We conclude that $\disc(\chi_A)$ is homogeneous of degree $n^2-n$ in 
the entries of $A$. 

We now apply Corollary~\ref{cor:hypersurface} 
with $p=n^2-1$ and $d=n^2-n$ to get (use (\ref{eq:wilki}))
\begin{equation}\tag*{\cuadrado}
   \bfE_{A\in B(M,\sigma)}(\ln \kappa_{\sf eigen}(A))
    \leq 8\ln n+2\ln\frac1\sigma+4 +\frac12\ln 2. 
\end{equation}

\subsection{Complex polynomial systems}

Let $d_1,\ldots,d_n\in\N\setminus\{0\}$. We denote by $\Hd$ the
vector space of polynomial systems $f=(f_1,\ldots,f_n)$ with
$f_i\in\C[X_0,\ldots,X_n]$ homogeneous of degree $d_i$,
$i=1,\ldots,n$. For $f,g\in \Hd$ we write
\begin{equation*}
   f_i(x) = \sum_{\alpha} a^i_\alpha X^\alpha, \quad g_i(x)
   = \sum_{\alpha} b^i_\alpha X^\alpha,
\end{equation*}
where $\alpha=(\alpha_0, \dots, \alpha_n)$ is assumed to range over
all multi-indices such that $|\alpha| = \sum_{k=0}^n \alpha_k = d_i$
and $X^\alpha:= X_0^{\alpha_0}X_1^{\alpha_1}\cdots X_n^{\alpha_n}$.

The space $\Hd$ is endowed with a Hermitian inner product $\langle
f,g \rangle = \sum_{i=1}^n \langle f_i, g_i \rangle$, where
\begin{equation*}
  \langle f_i, g_i \rangle =
  \sum_{|\alpha|=d_i}
   a^i_\alpha\, \bar{b^i_\alpha}\, {d_i \choose \alpha}^{-1}.
\end{equation*}
Here, the bar denotes complex conjugate and
the multinomial coefficients are defined by:
\begin{equation*}
  {d \choose \alpha} = \frac{d!}{\alpha_0 ! \alpha_1 ! \cdots \alpha_n !}.
\end{equation*}
Note that choosing this Hermitian product amounts to choosing the
monomials ${\tiny \sqrt{{d_i\choose \alpha}}}X^\alpha$ as
orthonormal basis of $\Hd$.

In the case of one variable, this product was introduced by H. Weyl
\cite{weyl:50}. Its use in computational mathematics goes back at
least to Kostlan~\cite{Kostlan93}. Throughout this
section, let $\|f\|$ denote the corresponding norm of $f$.
As described in \S\ref{ss:prelim}, the Weyl product defines a
Riemannian structure on the corresponding space $\proj(\Hd)$, with
its associated projective distance $d_{\proj(\Hd)}$.

In a seminal series of papers, M.~Shub and
S.~Smale~\cite{Bez1,Bez2,Bez3,Bez4,Bez5} studied the problem of,
given $f\in\Hd$, compute (an approximation of) a zero of $f$. They
proposed an algorithm and studied its complexity in terms of, among
other parameters, a condition number $\munorm(f)$ for $f$. We recall
its definition (see~\cite[Chapter~12]{bcss:95} for details).
For a simple zero $\zeta\in\proj^n$ of $f\in\Hd$ one defines
\begin{equation*}
  \munorm(f,\zeta):= \|f\|\left\| (Df(\zeta)_{|T_\zeta})^{-1}
  \mathsf{diag}(\sqrt{d_1}\|\zeta\|^{d_1-1},\dots,\sqrt{d_n}
   \|\zeta\|^{d_n-1})\right\|,
\end{equation*}
where $Df(\zeta)_{|T_\zeta}$ denotes restriction of the derivative of
$f\colon\C^{n+1}\to\C^n$ at $\zeta$ to the tangent space
$T_\zeta\proj^n =\{v\in\C^{n+1}\mid \langle v,\zeta\rangle = 0\}$
of $\proj^n$ at $\zeta$.
Note that $\munorm(f,\zeta)$ is homogeneous
of degree 0 in $f$ and $\zeta$.
If $f$ has only simple zeros $\zeta_1,\ldots,\zeta_q$ we define
\begin{equation*}
   \munorm(f):=\max_{i\leq q}\munorm(f,\zeta_i);
\end{equation*}
otherwise we set $\munorm(f)=\infty$.
The study of $\munorm(f)$ plays a central role in the series
of papers above. A main result is the
following~\cite{Bez2} (see
also~\cite[Theorem~1, Chapter~13]{bcss:95}).

\begin{theorem}\label{th:SS}
Let $n>1$. The probability that $\munorm(f)>1/\e$ for
$f\in\proj(\Hd)$ and $\e>0$ is less than or equal to
$$
   \e^4n^3(n+1)N(N-1)\DD
$$
where $\dim \Hd=N+1$ and $\DD=\displaystyle\prod_{i=1}^nd_i$
is the B\'ezout number.
\end{theorem}

We want to extend Theorem~\ref{th:SS} to a smoothed analysis
of $\munorm(f)$. To do so, we first bound $\munorm(f)$
by a conic condition number. Let $\Sigma\subset\proj(\Hd)$ be the
discriminant variety, which consists of the systems $f\in\proj(\Hd)$
having multiple zeros.
The Condition Number Theorem~\cite[\S12.4]{bcss:95} states that,
for a zero $\zeta\in\proj^n(\C)$ of $f$,
\begin{equation*}
   \munorm(f,\zeta)=\frac{1}{d_{\proj(\Hd)}(f,\Sigma\cap V_{\zeta})},
\end{equation*}
where $V_\zeta:=\{f\in\proj(\Hd)\mid f(\zeta)=0\}$. Therefore,
\begin{equation*}
   \munorm(f)=\max_{i\leq q}\munorm(f,\zeta_i)
   =\frac{1}{\min_{i\leq q} d_{\proj(\Hd)}(f,\Sigma\cap V_{\zeta_i})}
   \leq \frac{1}{d_{\proj(\Hd)}(f,\Sigma)}.
\end{equation*}
We can now proceed with the desired extension.

We identify the $f_i$ with their coefficient vectors
in $\C^{N_i}$, where $N_i= {n+d_i \choose d_i}$.
Set $N=\sum_i N_i-1$ so that $\Sigma\subset\proj^{N}$. 
Our next result bounds the degree of $\Sigma$. Similar 
bounds were given in~\cite[Proposition~6.1]{Ren87a}. 

\begin{lemma}\label{lem:disc}
The discriminant variety $\Sigma \subset\proj^{N}$ is a
hypersurface, defined by a multihomogeneous polynomial of 
total degree
$$
 \left(1 + \left(1-n + \sum_{i=1}^n d_i\right)
  \sum_{i=1}^n\frac{1}{d_i}\right)
  \DD \le 2n\DD^2
$$
in the coefficients of $f_1,\dots,f_n$.
\end{lemma}

\begin{Proof}
Given $n+1$ homogeneous polynomials $f_0,\dots,f_n$ in
$\C[X_0,X_1,\dots,X_n]$ of degrees $d_i$, it is known
(see~\cite[\S4.2]{stur:02})
that there exists an irreducible polynomial 
$\mathsf{res}(f_0,\dots,f_n)$
in the coefficients of the $f_i$ (unique up to a scalar)
such that $\mathsf{res}(f_0,\dots,f_n)=0$ if and only
if the system $f_0=\cdots f_{n}=0$ has a projective solution.
(The polynomial $\mathsf{res}$ is called the {\em multivariate resultant}.)
Moreover, $\mathsf{res}$ is multihomogeneous of degree
$\prod_{j\neq i}d_j$ in the coefficients of each $f_i$.

Now define
\begin{equation*}
  \d(f_1,\ldots,f_n):=\mathsf{res}(g,f_1,\dots,f_n),
\end{equation*}
where $g:=\det(df_1,\dots,df_n,\sum_i X_i dX_i)$.
A solution $\zeta$ to the system $f=0$ is degenerate if and only 
if the $df_i(\zeta)$ are linearly dependent, which is the case 
if and only if $g(\zeta)=0$ (here we used Euler's identity, 
stating that for homogeneous $f_i$ and all
$x\in\C^{n+1}$, $df_i(x)$ is orthogonal to $x$).
It follows that $\d(f_1,\dots,f_n)$ defines the discriminant 
variety $\Sigma$.

For the degree calculations, note first that
$\deg g=1+\sum_{i=1}^n (d_i-1) = 1-n + \sum_{i=1}^n d_i$.
We thus obtain
$$
 \deg\delta(f_1,\dots,f_n) 
 = \DD + \deg g \sum_{i=1}^n \frac{\DD}{d_i}
 = \DD \left( 1 + \deg g \sum_{i=1}^n d_i \right)
$$
as claimed. This degree can be (rather crudely) 
estimated by $2n\DD^2$.
\end{Proof}

\begin{theorem}
For all $f\in\proj(\Hd)$, all
$\sigma\in(0,1]$, and all $t\geq N\sqrt2$ we have
\begin{equation*}
   \Prob_{g\in B(f,\sigma}\{\munorm(g)\geq t\}\leq
    4 N^3 e^3 n \DD^2
   \left(\frac{1}{t\sigma}\right)^{2}
   \left(1+N\frac{1}{t\sigma}\right)^{2(N-1)}
\end{equation*}
and
\begin{equation*}
   \bfE_{g\in B(f,\sigma}(\ln\munorm(g))\leq
  \frac{7}{2}\ln N+\ln \DD +\frac12\ln n
   +5+2\ln\left(\frac1\sigma\right).
\end{equation*}
\end{theorem}

\begin{Proof}
It follows from Corollary~\ref{cor:hypersurface} and
Lemma~\ref{lem:disc}.
\end{Proof}

{\small
\bibliography{condition}
}
\end{document}